\documentclass[a4paper,oneside,11pt]{article} 

\usepackage{mystyle}

\title{A transportation approach to the mean-field approximation}  
\author{Fanny Augeri\footnote{Weizmann Institute of Science, Israel, E-mail: \href{mailto:fanny.augeri@weizmann.ac.il}{fanny.augeri@weizmann.ac.il.} This work was supported by the ERC advanced grant LogCorFields. } } 
\date{\today} 
\begin{document}
\maketitle

\abstract{We develop transportation-entropy inequalities which are saturated for measures such that their log-density with respect to the background measure is an affine function, in the setting of the uniform measure on the discrete hypercube and the exponential measure. In this sense, this extends the well-known result of Talagrand in the Gaussian case. By duality, these transportation-entropy inequalities imply a strong integrability inequality for Bernoulli and exponential processes.
As a result, we obtain a dimension-free mean-field approximation of the free energy of a Gibbs measure and a dimension-free nonlinear large deviations bound  on the discrete hypercube. Applied to the Ising model, we deduce that the mean-field approximation is within $O(\sqrt{n} ||J||_2)$ of the free energy, where $n$ is the number of spins and $||J||_2$ is the Hilbert-Schmidt norm of the interaction matrix.  Finally, we obtain a reverse log-Sobolev inequality on the discrete hypercube similar to the one proved recently in the Gaussian case by Eldan and Ledoux.
}

\section{Introduction}

A fundamental question in statistical Physics is to understand the behavior of Gibbs measures, in particular through the computation of their \textit{free energy}. If  $\mu$ is the uniform measure on the discrete hypercube $\{-1,1\}^n$ and  $f: \{-1,1\}^n \to \RR$ is a function, called the \textit{potential}, one can consider the \textit{Gibbs measure} associated to $f$, defined as the probability measure
$$ \nu = Z_f^{-1}e^{f}d\mu,$$
where $Z_{f} = \int e^{f } d\mu$ is the \textit{partition function} of $\nu$. The logarithm of the partition function is called the \textit{free energy}. The knowledge of the free energy for the family of  potentials $\beta f$ for $\beta>0$ encodes a rich information on the Gibbs measure. Unfortunately, the free energy is generally an intractable quantity, which in turn motivates the search for meaningful large $n$ approximations.  The \textit{Gibbs variational principle} (see \cite[Lemma 6.2.13]{DZ}) asserts that the free energy admits the following variational form
\begin{equation}\label{Gibbs} \log \int e^{f} d\mu = \sup_{\nu } \big\{ \int f d\nu - H(\nu|\mu) \big\},\end{equation}
where the supremum runs over all probability measures $\nu$ on $\{-1,1\}^n$, and $H(\nu|\mu)$ denotes the \textit{relative entropy} between $\nu$ and $\mu$. The \textit{mean-field approximation} consists in restricting the above supremum over the special class of product probability measures (or more generally \textit{tilted measures}, that is, measures whose log-density with respect to the background measure is an affine function). As product probability measures on the discrete hypercube are parametrized by their mean, the mean-field approximation reduces an optimization problem on probability measures on $\{-1,1\}^n$ into an optimization problem on $[-1,1]^n$, which is much more tractable. The question is then to understand under which condition on the potential $f$, the mean-field approximation can be justified rigorously. The Gibbs variational principle implies that the mean-field approximation always gives a lower bound on the free energy, that is
\begin{equation} \label{lowerbound} \log \int e^{f } d\mu \geq \sup_{ y \in [-1,1]^n } \big\{ \int f d \mu_y - I(y) \big\},\end{equation}
where  $\mu_y$ is the product measure on $\{-1,1\}^n$ with barycenter $y$, and $I(y) = H(\mu_y | \mu)$.
Another way to reformulate the accuracy of the mean-field approximation is to say that the above inequality is approximately tight in the large $n$ limit. Our main task in the present work will be to obtain quantitative upper bounds.

In a seminal paper \cite{CD}, Chatterjee and Dembo showed that given an extension of the potential $f$ to the hypercube $[-1,1]^n$, the mean-field approximation is accurate if the set of gradients of $f$ is of low complexity in a $\ell^2$-metric entropy sense. However, the quantitative error bound from the mean-field approximation they obtained is rather intricate, and involves in particular $L^{\infty}$-norms of the partial derivatives of $f$ up to the second order.  

In the case of the Ising model, where the potential $f$ is a quadratic form
$$\forall x \in \{-1,1\}^n, \  f(x) = \langle x, J x \rangle,$$
given in terms of an \textit{interaction matrix} $J$, the convergence of the free energy to the mean-field approximation was shown in the context of dense graphs using the graphon framework in \cite{BCLSV} and \cite{BCCZ}. For general graphs,
a first breakthrough was made by Basak and Mukherjee \cite{BM} who showed the accuracy of the mean-field approximation under the condition that $ ||J||_2 = o(\sqrt{n})$, denoting by $||J||_2$ the Hilbert-Schmidt norm of $J$. However, this result does not give any information about the speed of convergence. A quantitative error bound from the mean-field approximation in $O( (n||J||_2)^{\frac{2}{3}})$ up to a logarithmic factor was derived by Jain, Koehler and Mossel \cite{JKM} using the Frieze-Kannan regularity lemma.

Another approach to this problem goes through the decomposition of the Gibbs measure itself into a mixture of measures where the coordinates are weakly correlated. This line of research was exploited by Jain, Koehler and Risteski in \cite{JKR} to remove the logarithmic correction in the mean-field approximation for the Ising model, and showed that,
\begin{equation} \label{JKR} \log \int e^{\langle x, Jx\rangle } d\mu(x) \leq \sup_{y\in[-1,1]^n } \{ \langle y ,Jy \rangle - I(y)\} + O\big( (n ||J||_2)^{\frac{2}{3}}\big).\end{equation}
In \cite{Eldan}, Eldan proved a structural theorem for general Gibbs measures in Gaussian space and for the discrete hypercube. He deduced an upper bound on the free energy where the complexity of the discret gradient of the potential is assessed in terms of its \textit{Gaussian mean-width}, namely,
\begin{equation} \label{Gaussianwidth} g(V) = \EE \sup_{\xi \in V} \langle \xi, \Gamma \rangle,\end{equation}
where $V = \nabla f(\{-1,1\}^n)$ and $\Gamma$ is a standard Gaussian variable in $\RR^n$. His approximation of the free energy \cite[Corollary 2]{Eldan} takes the form, 
\begin{equation} \label{Eldanmeanfield} \log \int e^f d\mu \leq \sup_{ y\in [-1,1]^n} \big \{ \int f d\mu_y - I(y)\big\} + O\big( \mathrm{Lip}(f)^{\frac{2}{3}} g(V)^{\frac{1}{3}} n^{\frac{2}{3}} \big),\end{equation} 
where $\mathrm{Lip}(f)$ is the Lipschitz constant of $f$ with respect to the Hamming metric. This approach was further developed by Austin \cite{A} who proved a structure theorem for Gibbs measures on general product spaces and deduced a mean-field approximation of the free energy.

In \cite{NL}, the author proved a mean-field approximation for Gibbs measures with respect to general compactly supported background measures which, using Sudakov minoration, implies in the case of the discrete hypercube that
$$ \log \int e^{f} d\mu \leq \sup_{ y\in [-1,1]^n} \big\{ \int fd \mu_y - I(y)\big\}  + O\big( g(V)^{\frac{2}{3}} n^{\frac{1}{3}}\big).$$
In particular, this bound enables one to recover the bound \eqref{JKR} for the Ising model. In the present paper, we will remove the dimension dependence from the above estimate, and prove the dimension-free inequality,
\begin{equation}\label{dimfree} \log \int e^{f} d\mu \leq \sup_{ y\in [-1,1]^n} \big\{ \int fd \mu_y - I(y)\big\}  + O\big(b(V)\big),\end{equation}
where $b(V) = \EE \sup_{\xi \in V} \langle \xi, \eps\rangle$, and $\eps$ is uniformly distributed on $\{-1,1\}^n$.

In a recent work \cite{Eldanmean}, Eldan proved a new decomposition theorem which allowed him to show in the case of the Ising model that for any $p>0$, the error on the free energy induced by the mean-field approximation is $O( \frac{1+p}{p} (n||J||_p)^{\frac{p}{1+p}})$, where $|| \ ||_p$ denotes the $p$-Schatten norm. This bound recovers for $p=2$ the previous $O((n ||J||_2)^{\frac{2}{3}})$ error shown by  Jain, Koehler and Risteski in \cite{JKR}, and can significantly improve upon this bound by an appropriate choice of $p$.

The goal of this paper is to propose a transportation approach for the problem of the mean-field approximation of the free energy of Gibbs measures in the specific case of the discrete hypercube. The main interest of this approach is that it allows us to derive an approximation which is dimension-free, i.e \eqref{dimfree}. We develop new transportation-entropy inequalities in the case of the Bernoulli and the exponential distribution. 
Originally, the transportation-entropy inequalities were put forward by Marton \cite{Marton} and Talagrand \cite{Talagrand2}. They appear to have strong connections with concentration inequalities (see  \cite[Chapter 6]{Ledouxmono}, \cite{Dembo}, \cite[Chapter 8]{BLM}, or \cite[section 4]{GL}). They also have many links with other functional inequalities. Quadratic transportation-entropy inequalities are known to imply a spectral gap inequality by \cite[section 4.1]{BGL} (see also \cite[section 8.3]{GL}), and are weaker that logarithmic Sobolev inequalities by the result of Otto and Villani \cite{OV} (see also \cite{BGL}).

The main feature of the transportation-entropy inequalities we will present is that they are saturated by \textit{tilts} of the background measure, that is  measures with an affine log-density. Given the central role of such probability measures in the mean-field approximation, this feature will be particularly crucial. 

Using duality, the main consequence we derive from these transportation-entropy inequalities is a strong integrability inequality for Bernoulli and exponential processes, similar to the Gaussian case. In turn, this will provide us the main ingredient to obtain a dimension-free mean-field approximation of the free energy of Gibbs measures and in a similar fashion, a dimension-free nonlinear large deviations bound on the discrete hypercube.
In the setting of the Ising model on $\{-1,1\}^n$, we deduce that the mean-field approximation is within $O( \sqrt{n} ||J||_2)$ of the free energy, improving  the previous known bound \eqref{JKR} involving the Hilbert-Schmidt norm of $J$. Finally, we prove a dimension-free reverse log-Sobolev inequality on the discrete hypercube similar as the one existing in the Gaussian case \cite{LedouxE}.


\section{Main results}

\subsection{Transportation-entropy inequalities}
Let $\mathcal{P}(\RR^n)$ denote the set of probability measures on $\RR^n$. For any $\mu, \nu \in \mathcal{P}(\RR^n)$, and a lower semi-continuous cost function $c : \RR^n \times \RR^n \to [0,+\infty]$, one defines the \textit{transportation cost} $\mathcal{W}_c(\nu,\mu)$ by,
$$ \mathcal{W}_c(\nu,\mu)= \inf_{\pi} \int c(x,y) d\pi(x,y),$$
where the infimum runs over all couplings between $\nu$ and $\mu$. We say that a given measure $\mu$ satisfies a \textit{transportation-entropy inequality} with cost function $c : \RR^n \to [0,+\infty]$ if,
\begin{equation} \label{transpoin} \forall \nu \in \mathcal{P}(\RR^n), \   \mathcal{W}_c(\nu,\mu) \leq H(\nu|\mu),\end{equation}
where $H(\nu|\mu)$ denotes the relative entropy. 

Let $\mu$ be a reference  probability measure on $\RR^n$. We call $\nu \in \mathcal{P}(\RR^n)$ a \textit{tilt} of $\mu$ if the log-density with respect to $\mu$,  $\log \frac{d \nu}{d \mu}$, is an affine function. We address the question of finding a transportation-entropy inequality which is saturated by tilts of the reference measure $\mu$.  By Talagrand's result \cite{Talagrand2}, we know that the standard Gaussian measure on $\RR^n$, which we denote by $\gamma$, satisfies a transportation-entropy inequality with cost function $(x,y)\mapsto \frac{1}{2} ||x-y||_{\ell^2}^2$, where $|| \ ||_{\ell^2}$ denotes the $\ell^2$-norm, that is,
\begin{equation} \label{TalG} \forall \nu \in \mathcal{P}(\RR^n), \ \mathcal{W}_{\frac{1}{2} || \ ||_{\ell^2}^2 } (\nu,\gamma) \leq H(\nu|\gamma).\end{equation}
 As one can observe, this transportation-entropy inequality is tight for \textit{tilts} of the Gaussian measure, which are just push-forwards by translations.

In the case of the exponential measure $\eta =e^{-x}\Car_{x\geq 0} dx$ , we consider the following cost function,
\begin{equation} \label{defcoutexp}  \forall x,y \in \RR^n, \ c(x,y) = \sum_{i=1}^n  y_i \Lambda^*_{\eta}\Big( \frac{x_i}{y_i}\Big),\end{equation}
where for any $t >0$,
$$\Lambda^*_\eta (t) = t- 1-\log t,$$
and for $t\leq 0$, $\Lambda^*_\eta(t) = +\infty$.
With these definitions, we have the following transportation-entropy inequality.

\begin{Pro}\label{transpoexp}
Let $\eta$ be the probability measure $e^{-x}\Car_{x\geq 0} dx$ and $\eta^n$ be its $n$-fold product. For any probability measure $\nu$ on $\RR^n$,
$$\mathcal{W}_c(\nu,\eta^n) \leq H(\nu| \eta^n).$$
Moreover, the equality holds if $\nu$ is a tilt of $\eta^n$.
\end{Pro}
\begin{Rem}
In \cite{Talagrand2}, Talagrand proved that the symmetric exponential measure $m = \frac{1}{2}e^{-|x| } dx$ satisfies a certain family of transportation-cost inequalities with costs $c_t$ indexed by $t \in (0,1)$, defined by,
$$\forall x \in \RR, \  c_t(x,y) = \Big(\frac{1}{t}-1\Big) \big( e^{-t |x-y|} + t |x-y| -1\big).$$
This family of cost functions has the striking property that for \textit{any} $t$ there exists a probability measure which achieves the inequality in the transportation-cost inequality with cost function $c_t$. In this sense, this is a family of optimal cost functions. However, the probability measures which saturate the inequality are not tilts of the exponential measure, but are more intricate measures whose monotonous rearrangements from the exponential measure satisfy a certain family of differential equations.

\end{Rem}

To deal with the singularity of the Bernoulli measure, we propose a variant of the transportation-entropy inequality \eqref{transpoin} where we make it possible to enrich the transportation problem by considering another measure than the reference measure. More precisely, we will say that $\mu$ satisfies a transportation-entropy inequality if 
$$  \forall \nu \in \mathcal{P}(\RR^n), \   \mathcal{W}_c(\nu,\tilde \mu) \leq H(\nu|\mu),$$
where $c : \RR^n \times \RR^n \to [0,+\infty]$ is lower semi-continuous and $\tilde \mu$ is a fixed probability measure on $\RR^n$.

Let $p\in(0,1)$ and $I_p$ be the function defined by
\begin{equation} \label{defIp}\forall x \in[-1,1], \ I_p(x) =  \frac{1+x}{2} \log \frac{1+x}{2p} + \frac{1-x}{2} \log \frac{1-x}{2(1-p)},\end{equation}
and $I_p(x) = +\infty$ otherwise.
We define the cost function $w_p : \{-1,1\}^n \times [-1,1]^n\to [0+\infty]$ by 
\begin{equation} \label{defw} \forall  x \in \{-1,1\}^n, u \in [-1,1]^n, \ w_p(x,u) = \sum_{i=1}^n 2 |I'_p(u_i)| \Car_{x_i (h_0- u_i) <0},\end{equation}
and $h_0 = 2p-1$. 
With these definitions, we have the following transportation-entropy inequality.
\begin{Pro}\label{transpo} Let $\mu_p= (1-p) \delta_{-1} + p \delta_{1}$ and $\mathcal{U}$ be the uniform probability measure on $[-1,1]$.
For any probability measure $\nu$ on $\{-1,1\}^n$,
$$ \mathcal{W}_{w_p}( \nu,\mathcal{U}^n) \leq H(\nu | \mu_p^n),$$
and equality holds if $\nu$ is a product measure.
\end{Pro}
\subsection{Strong integrability of empirical processes}
The first consequence we will derive from the transportation-entropy of the previous section consists in the strong integrability of Bernoulli and exponential empirical processes. By empirical process, we mean any process of the form,
$$ (\langle \xi, X \rangle)_{\xi \in V},$$
where $V$ is some countable subset of $\RR^n$, $X$ is a random vector in $\RR^n$ with independent and identically distributed coordinates, and $\langle . , . \rangle$ denotes the standard inner product in $\RR^n$.

In the Gaussian case, it is known that for any countable set $V \subset \RR^n$,
\begin{equation} \label{strongG} \log \int e^{\sup_{\xi \in V} \{ \langle \xi, x \rangle - \frac{1}{2} ||\xi||_{\ell^2}^2\} } d \gamma(x) \leq \int \sup_{\xi \in V} \langle \xi, x \rangle d\gamma(x).\end{equation}
This inequality was first put forward in \cite{Vitale}. It can also be seen as a consequence of Talagrand's transportation-entropy for the Gaussian measure \eqref{TalG}.

We show that a similar estimate holds for the uniform measure on $\{-1,1\}^n$ and the exponential measure, where the quadratic cost $\frac{1}{2} ||\xi||_{\ell^2}^2$ is replaced by the \textit{logarithmic Laplace transform} of the measure considered.
In the following, for any probability measure $\mu$ on $\RR^n$ we will denote by $\Lambda_\mu$  its \textit{logarithmic Laplace transform}, defined as
$$ \forall \xi \in \RR^n , \ \Lambda_{\mu}(\xi) = \log \int e^{\langle \xi, x \rangle} d\mu(x).$$
For a subset $V \subset \RR^n$, we will also denote by $b(V)$ the \textit{``Rademacher mean-width''} of $V$, defined by
\begin{equation} \label{Rademwidth} b(V) = \EE \sup_{\xi \in V} \langle \xi, \eps \rangle,\end{equation}
where $\eps$ is uniformly sampled on $\{-1,1\}^n$. With this notation, we have the following result in the discrete setting. This estimate will be the key element of our approach to the mean-field approximation.
\begin{Pro}\label{strongintB}
Let $\mu$ be the uniform probability measure on $\{-1,1\}^n$. There exists a universal constant $\kappa >0$, such that for any $V \subset \RR^n$, 
$$ \log \int e^{\sup_{\xi \in V}\{ \langle \xi, x \rangle - \Lambda_{\mu}(\xi) \}}d\mu(x) \leq \kappa b(V).$$
\end{Pro}
Similarly, we get in the case of the exponential measure the following result. 
\begin{Pro}\label{intexpo}
Let $\eta$ be the probability measure $\Car_{x\geq 0} e^{-x} dx$. For any countable subset $V \subset \RR^n$, 
$$ \log \int e^{\sup_{\xi \in V} \{ \langle \xi, x\rangle - \Lambda_{\eta^n}(\xi) \}}d\eta^n(x) \leq \int \sup_{\xi \in V} \langle\Lambda_\eta(\xi), x-u\rangle d\eta^n(x),$$
where $u$ denotes the vector $(1,1,\ldots,1)$, and $\Lambda_\eta(\xi ) =(\Lambda_\eta(\xi_1),\ldots,\Lambda_\eta(\xi_n))$.
\end{Pro}

\subsection{Mean-field approximation}
Building on the previous strong integrability inequality for Bernoulli empirical processes, we prove a dimension-free mean-field approximation of the free energy of Gibbs measures. In the following we denote by $I$ the function defined by,
\begin{equation} \label{defI} \forall x \in [-1,1]^n, \   I(x) =\sum_{i=1}^n \Big( \frac{1+x_i}{2} \log (1+x_i) + \frac{1-x_i}{2} \log (1-x_i)\Big),\end{equation}
and $I(x)= +\infty$ otherwise. With this notation, we have the following result.
\begin{The}\label{meanfield} Let $\mu$ be the uniform probability measure on $\{-1,1\}^n$. There exists a universal constant $\kappa>0$, such that 
for any function $f : \RR^n \to \RR$  continuously differentiable,
$$ \log \int e^f d \mu \leq \sup_{ y\in [-1,1]^n } \{ f(y) - I(y) \} + \kappa b(V),$$
where $V = \nabla f([-1,1]^n)$, and $b(V)$ is defined in \eqref{Rademwidth}.
\end{The}
 One can interpret $I$ by the identity
$$\forall y \in [-1,1]^n, \  I(y) = H(\mu_y|\mu),$$
where $\mu_y$ stands for the unique product probability measure on $\{-1,1\}^n$ with mean $y$.
\begin{Rem}\label{harmoext}
If $f : \{-1,1\}^n \to \RR$ is some function defined only on the discrete hypercube, one can extend it harmonically to $[-1,1]^n$ by the formula
$$\forall y \in [-1,1]^n, \  f(y) = \int f d\mu_y.$$
We know by \cite[Fact 14]{Eldan} that
 $$ \nabla f(h) = \int \nabla f d \mu_y,$$
where the gradient on the right-hand side is the discrete gradient of $f$, that is,
$$\forall x \in \{-1,1\}^n,\    \nabla f(x) = (\partial_1 f(x),\ldots,\partial_n f(x)),$$
with 
$$\forall i \in \{1,\ldots,n\}, \ \partial_i f(x) = \frac{1}{2} \big( f(x_+) - f(x_-)\big),$$
where $x_+ = (x_1,\ldots ,x_{i-1},1,x_{i+1},\ldots ,x_n)$ and $x_- = (x_1,\ldots ,x_{i-1},-1,x_{i+1},\ldots,x_n)$.  

Thus, for this extension, the set of gradients 
$\nabla f([-1,1]^n)$ is the convex hull of the discrete gradients. Therefore, the error term of Theorem \ref{meanfield} is just $b(\nabla f(\{-1,1\}^n)$).

\end{Rem}

\begin{Rem}
If $f :\RR^n \to \RR$ is a continuously differentiable function and $g$ is the harmonic extension of $f_{|\{-1,1\}^n}$ to $[-1,1]^n$ described in the previous remark, then 
$$ \sup_{y \in[-1,1]^n} \big( f(y) - g(y)\big)  \leq \lambda  b(V),$$
where $V = \nabla f( [-1,1]^n)$, and $\lambda$ is a numerical constant. We refer the reader to Lemma \ref{exten} for a proof of this fact.
This implies that in the mean-field approximation stated in Theorem \ref{meanfield}, the extension one chooses does not matter as soon as its set of gradient is of low complexity. 

\end{Rem}
Applying Theorem \ref{meanfield} to the Ising model, we obtain the following corollary.
\begin{Cor}\label{Isingmean}
 Let $J$ be a Hermitian matrix of size $n$ such that $J_{i,i} = 0$ for any $i\in\{1,\ldots,n\}$, and $h\in \RR^n$. Then,
$$ \log \int e^{\langle x, Jx \rangle+ \langle h, x\rangle }d\mu(x)  \leq \sup_{y \in [-1,1]^n} \{ \langle y,J y\rangle + \langle h,y \rangle  - I(y)\} +  \kappa || J ||_2 \sqrt{n},$$
where $\kappa$ is a universal positive constant, and $|| \ ||_2$ denotes the Hilbert-Schmidt norm, namely
$$ || J||_2 = \Big( \sum_{1\leq i,j \leq n} |J_{i,j}|^2 \Big)^{1/2}.$$
\end{Cor}

In the interesting large deviation regime, the free energy is expected to be of order $n$. Thus, the above Corollary  \ref{Isingmean} gives a meaningful upper bound whenever $||J||_2 = o(\sqrt{n})$. It recovers the qualitative result of Basak and Mukherjee \cite{BM} for the Ising model, and gives a quantitative error term which is strictly smaller than the one found in \cite{JKR}, i.e \eqref{JKR}.

%
%

\begin{Exe}[$d$-regular graphs] Let us consider a $d$-regular graph $G$ with $n$ vertices ($n$ and $d$ are implicitly taken such that $nd$ is even). Let $A$ denote the adjacency matrix of $G$, and let us consider the Ising model with the interaction matrix $J = \frac{1}{d} A$. This scaling is taken so that the free energy of this model is of order $n$. 
As $||J||_2 = \sqrt{n/d}$, Corollary \ref{Isingmean} gives
$$ \log \int e^{\langle x, Jx \rangle} d\mu(x) = \sup_{y \in [-1,1]^n} \{ \langle y,J y\rangle -I(y)\} +  O\Big(\frac{n}{\sqrt{d}}\Big).$$
This bound improves the one of Eldan \cite[Example 3]{Eldanmean} who showed that for a $d$-regular expander, that is such that the second largest eigenvalue $\lambda_2(A) =O(\sqrt{d})$, one has the error 
$$ \frac{n}{\sqrt{d}^{1-o(1)}},$$
over the mean-field approximation, in the regime $\log d \ll \log n$.

\end{Exe}

\begin{Rem}
In \cite{Eldanmean}, Eldan proved that for any $p>0$ the mean-field approximation is within  $O( \frac{p}{p+1} (n||J||_p)^{\frac{p}{p+1}})$ of the free energy of the Ising model with interaction matrix $J$, where $||\ ||_p$ is the $p$-Schatten norm.
One can note that for $p\geq 2$ and in the regime where Eldan's bound is meaningful, that is $||J||_p =O( n^{1/p})$, the inequality $|| J||_2 \leq n^{\frac{1}{2}-\frac{1}{p}} ||J||_p$, yields that the error term of Corollary \ref{Isingmean} is smaller than $O( (n||J||_p)^{\frac{p}{p+1}})$. However, the real interest of Eldan's bound is when $p\leq 2$ and in particular the regime when $p\to 0$ when $n$ grows to infinity. For $p<2$, it seems that Eldan's bound and the one given by Corollary \ref{Isingmean} cannot be compared in general. For specific examples, like the Curie-Weiss model or the lattice with mesoscopic interactions (see \cite[example 2]{Eldanmean}) where the eigenvalues of the interaction decrease exponentially fast to $0$, Eldan's bound is  better and yields only logarithmic errors, whereas Corollary \ref{Isingmean} can only provide error terms depending polynomially in the dimension. 
\end{Rem}

\subsection{Nonlinear large deviations} 
The theory of nonlinear large deviations was introduced by Chatterjee and Dembo \cite{CD} in order to understand the large deviations of nonlinear functions of independent Bernoulli random variables. One of the motivation for developing this theory comes from the question of the deviations of sub-graph counts in sparse Erd\H{o}s--R\'enyi graphs. 

Given a function $f : \{-1,1\}^n \to \RR$, and $X$ uniformly sampled on $\{-1,1\}^n$, one can wonder when the optimal change of measure in the large deviations of $f(X)$ is given by product measures. The nonlinear large deviations theory aims at answering this question and at identifying which condition on $f$ can guarantee this mechanism of deviation to happen. Similarly as for the question of the mean-field approximation of the free energy, Chatterjee and Dembo showed in \cite{CD} that a sufficient condition is that the set of gradients of $f$ is of low complexity in a $\ell^2$-metric entropy sense. 

Efforts have been put into improving the original non-asymptotic bound of \cite{CD}, which has the inconvenient of involving error terms related to the smoothness of $f$. In \cite{Yan}, Yan generalizes to products of general compact spaces the nonlinear large deviations bound of \cite{CD}. Eldan \cite{Eldan} removed most of the smoothness assumptions and provided a bound where the complexity of the gradient is assessed in term of its Gaussian mean-width. In \cite[Corollary 2.2]{CoD}, Cook and Dembo proposed a nonlinear large deviation bound which has the specificity of not relying on the complexity of the gradient  but rather on an efficient covering of  the space by convex sets.

We propose here a nonlinear large deviations bound in the specific case of the discrete hypercube whose main feature is to be dimension-free. As we will show, it follows from the strong integrability inequality of Bernoulli processes of Proposition \ref{strongintB}. 

To describe this bound, we extend $I_p$, defined on $\RR$ by the formula \eqref{defIp}, to $\RR^n$ by setting,
$$ \forall y \in \RR^n, \ I_p(y) := \sum_{i=1}^n I_p(y_i).$$
Let  $f : [-1,1]^n \to \RR$ be a function, and define the rate function
$$ \forall t \in \RR, \ \phi_p(t) = \inf\{ I_p(y) : f(y) \geq t, y \in \RR^n \}.$$
With this notation, we have the following theorem.
\begin{The}\label{NL}Let $t \in \RR$ and $\delta>0$. Assume that
$$ \forall s>t-\delta, \ \phi_p(s) > \phi_p(t-\delta).$$
 Let $V = \nabla f([-1,1]^n)$ and let $X$ be a random vector sampled according to $\mu_p^n$. There exist universal constants $C,\kappa>0$, such that if
$$ b(V) \leq  \delta/\kappa ,$$
where $b(V)$ is defined in \eqref{Rademwidth}, then 
$$ \log \PP\big( f(X) \geq t \big) \leq - \phi_p(t-\delta) +C\log \Big(\frac{n L \log(1/p(1-p))}{\delta}\Big),$$
where $L= \sup_{x\in[-1,1]^n} || \nabla f(x)||_{\ell^2}$.
\end{The}

\begin{Rem}
It is actually possible to weaken the regularity assumption on $f$, and assume that for any $x,y \in [-1,1]^n$,
$$ f(x) - f(y) \leq \sup_{\xi\in V} \langle \xi, x-y\rangle,$$
where $V$ is a convex subset of $\RR^n$.
\end{Rem}

\subsection{Reverse log-Sobolev inequality on the discrete hypercube}
Let $\mu$ be the uniform measure on $\{-1,1\}^n$. The \textit{logarithmic Sobolev inequality}  on the discrete hypercube (see \cite[Theorem 5.1]{BLM}) says that for any $\nu = e^f d\mu$ probability measure on $\{-1,1\}^n$,
\begin{equation}\label{logSob}  H(\nu| \mu) \leq \frac{1}{2} \int || \nabla f (x) ||_{\ell^2}^2 d\nu(x),\end{equation}
where $\nabla f$ denotes the discrete gradient. The inequality \eqref{logSob} can be improved by replacing the quadratic function $|| \ ||_{\ell^2}^2/2$ by $I(\nabla \Lambda_{\mu})$, which gives
\begin{equation} \label{logSob2} H(\nu| \mu) \leq  \int I( \nabla \Lambda_{\mu} ( \nabla f (x))) d\nu(x).\end{equation}
From the inequality in dimension $1$,
$$ \forall \lambda \in \RR, \ I(\Lambda_{\mu}'(\lambda)) \leq \frac{1}{2} \lambda^2,$$
we see that  \eqref{logSob} is indeed implied by \eqref{logSob2}. 

The proof of inequality \eqref{logSob2} goes over induction on the dimension. For $n=1$, it is straightforward to see that there is actually equality. For $n>1$, one uses the sub-additivity of the relative entropy \cite[Proposition 5.6]{Ledouxmono},
$$ H(\nu|\mu) \leq \sum_{i=1}^n \int H(\nu_{x^{(i)}} | \mu_{1/2}) d\nu(x),$$
where $\nu_{x^{(i)}}$ is the conditional probability measure given $x^{(i)} = (x_1,\ldots,x_{i-1},x_{i+1},\ldots,x_n)$, which is equal in our case to the probability measure proportional to $e^{\partial_if(x)x_i} d \mu_{1/2}(x_i)$

The interest of stating the log-Sobolev inequality this way is that it is saturated for product measures. Thus, one can expect that whenever the gradient of $f$ is of low complexity, the inequality \eqref{logSob2} above is almost an equality. We will prove that it is indeed the case, and show a reverse log-Sobolev inequality, similar to the one proved by Eldan and Ledoux \cite{LedouxE} in the Gaussian case.

\begin{Pro}
\label{revlogSob}
Let $\nu =e^f d\mu$ be a probability measure on $\{-1,1\}^n$. Let
$$ \mathcal{I}(\nu) = \int I(\nabla \Lambda_\mu (\nabla f(x))) d\nu(x),$$
where $\nabla f$ is the discrete gradient of $f$.
 Then,
$$ \mathcal{I}(\nu) \leq H(\nu |\mu) + \kappa \int \sup
_{y \in \mathcal{C}_n} \langle \nabla f(y), x \rangle d\mu(x),$$ 
where $\kappa$ is a universal constant.

\end{Pro}

\section{Transportation-entropy inequalities} 
In this section, we prove the transportation-entropy inequalities of the Propositions \ref{transpoexp} and \ref{transpo} which are at the base of our results. We start by recalling some standard features of these inequalities.
An important property is that they tensorize in a certain way which we recall in the following lemma.  The reader may find a proof of this result in \cite[Proposition 1.3]{GL}.

\begin{Lem} \label{tensor}
If for $i\in \{1,2\}$,  $\mu_i, \tilde \mu_i$ are probability measures on $\RR^{d_i}$, such that
$$ \forall \nu \in \mathcal{P}(\RR^{d_i}), \ \mathcal{W}_{c_i}(\nu, \tilde \mu_i) \leq H(\mu| \mu_i),$$
where $c_i : \RR^{d_i} \times \RR^{d_i}\to [0,+\infty]$ is a lower semi-continuous function, then 
$$ \forall \nu \in \mathcal{P}(\RR^{d_1}\times \RR^{d_2}), \ \mathcal{W}_{c_1 \oplus c_2}(\nu, \tilde \mu_1\otimes \tilde \mu_2) \leq H(\nu | \mu_1\otimes \mu_2),$$
where $c_1 \oplus c_2$ is defined for any $x=(x_1,x_2) \in \RR^{d_1} \times \RR^{d_2}$ and $y =(y_1,y_2) \in \RR^{d_1} \times \RR^{d_2}$ by,
$$ c_1 \oplus c_2(x,y) = c(x_1,y_1) + c(x_2,y_2).$$
\end{Lem}
Note that when $\mu$ is a product measure, the tilts of $\mu$ are also product measures. Therefore, the question of finding a transportation-entropy inequality which is saturated for tilts reduces itself to a $1$-dimensional problem by the tensorization property described above together with the following fact. It follows from the definitions of the relative entropy and the transportation cost.

\begin{Fact}\label{tensoreq}
For $i\in \{1,2\}$, let $\mu_i$ and $\nu_i$  be probability measures on $\RR^{d_i}$, and $c_i : \RR^{d_i} \times \RR^{d_i} \to[0,+\infty]$ be lower semi-continuous functions. Then,
$$ H(\nu_1 \otimes \nu_2 | \mu_1 \otimes \mu_2) = H(\nu_1|\mu_1) + H(\nu_2|\mu_2),$$
and
$$\mathcal{W}_{c_1 \oplus c_2} (\nu_1 \otimes \nu_2, \mu_1 \otimes \mu_2) = \mathcal{W}_{c_1}(\nu_1,\mu_1) + \mathcal{W}_{c_2}(\nu_2,\mu_2).$$

\end{Fact}

The main aspect of transportation-entropy inequalities we will use is their dual functional form, which consists of \textit{infimum-convolution inequalities}. This duality relies on the Kantorvitch duality \cite[Theorem 5.10]{Villani} which states that if $c : \RR^n \times \RR^n \to [0,+\infty]$ is lower semi-continuous, then for any $\nu,\mu \in \mathcal{P}(\RR^n)$,
\begin{equation} \label{Kantduality} \mathcal{W}_c(\nu, \mu) = \sup_{ \phi \in L^1(\nu)} \big\{ \int \phi d \nu - \int \phi^c d\mu\big\},\end{equation}
where $\phi^c$ the $c$-conjugate of $\phi$ defined by,
$$\forall y\in \RR^n,\ \phi^c(y) = \sup_{ x\in \RR^n} \{ \phi(x) -c(x,y)\}.$$
Moreover, by \cite[Theorem 5.10, (ii)]{Villani}, if $\mathcal{W}_c(\nu,\mu) <+\infty$, then a coupling $\pi$ between $\nu$ and $\mu$ is optimal if and only if there exists $\phi \in L^1(\nu)$ such that $\pi$-almost surely,
\begin{equation} \label{Kantdualityopt} \phi(y) - \phi^c(x) = c(x,y).\end{equation}

 In the next lemma, we recall the equivalence between transportation-entropy and infimum-convolution inequalities (see \cite[Corollary 3.1]{GL} or \cite[Theorem 5.26]{Villani}) and we show that an equality case in the transportation-entropy inequality can be translated into an equality case for the infimum-convolution inequality.

\begin{Pro}\label{duality}
Let $\mu,\tilde \mu \in \mathcal{P}(\RR^n)$ and $c : \RR^n \times \RR^n \to [0,+\infty]$ be a lower semi-continuous function. The following statements are equivalent.

\begin{enumerate}
\item[(i).] $\mu$ satisfies the transportation-entropy inequality,
\begin{equation} \label{transpoent}\forall \nu \in \mathcal{P}(\RR^n), \  \mathcal{W}_c(\nu,\tilde \mu) \leq H(\nu| \mu). \end{equation}

\item[ (ii).] $\mu$ satisfies the infimum-convolution inequality,
\begin{equation} \label{infconv} \log \int e^f d\mu \leq \int  f^c d \tilde \mu,\end{equation}
for any $f : \RR^n \to \RR$ measurable such that $f^c \in L^1(\mu)$.
\end{enumerate}
Let $\nu$ be such that $H(\nu|\mu)<+\infty$. In the case $(i)$ or $(ii)$ is satisfied, equality holds in \eqref{transpoent} for $\nu$ if and only if  equality holds in \eqref{infconv} for $f = \log \frac{d \nu}{d\mu}$.
\end{Pro}

\begin{proof}
A proof of the equivalence between $(i)$ and $(ii)$ can be found in  \cite[Corollary 3.1]{GL} or \cite[Theorem 5.26]{Villani}. We are now left to prove the equivalence between the equality cases. Assume $(i)$ holds and there is equality in \eqref{transpoent} for $\nu$. Denote by $f = \log \frac{d \nu}{d\mu}$. As $H(\nu|\mu) < +\infty$, we have that $\mathcal{W}_c(\nu,\tilde \mu) < +\infty$. By Kantorovich duality \eqref{Kantdualityopt}, there exists $\phi \in L^1(\nu)$ such that, $\pi$-almost surely,
\begin{equation}\label{Kant} \phi(x) - \phi^c(y) = c(x,y).\end{equation}
Our goal is to show that $\phi$ is equal to $f$ up to some additive constant. Indeed, if $\phi = f +\alpha$ for some constant $\alpha\in \RR$, then using the fact  that $(f+\alpha)^c = f^c +\alpha$, we have $\pi$-almost surely,
$$ f(x) -f^c(y) = c(x,y).$$
Integrating the above equality with respect to $\pi$ yields,
$$ \int f^c(y) d\tilde \mu(y) = \int f(x) d \nu(x) - \int c(x,y) d\pi(x,y).$$
As $\nu$ achieves the equality in \eqref{transpoent}, we obtain
$$ \int f^c(y) d\mu(y) = 0,$$
which would prove the first part of the equivalence between the equality cases.

Using Kantorovich duality \eqref{Kantduality}, we have for any $\eta \in \mathcal{P}(\RR^n)$ such that $\phi \in L^1(\eta)$,
$$ \mathcal{W}_c(\eta,\tilde \mu) \geq \int \phi d\eta - \int \phi^c d\tilde \mu.$$
Integrating  \eqref{Kant} with respect to $\pi$ and combining with the above inequality, we deduce
$$ \mathcal{W}_c(\eta,\tilde \mu) - \mathcal{W}_c(\nu,\tilde \mu) \geq \int \phi d(\eta - \nu).$$
But, as $(i)$ holds and equality holds for $\nu$ in \eqref{transpoent}, we get
\begin{equation}\label{subdiff} H(\eta|\mu) - H(\nu|\mu) \geq  \int \phi d(\eta - \nu).\end{equation}
Let $\psi : \RR^n \to\RR$ be a measurable and bounded function such that 
$$ \int \psi d\nu =0.$$
For $\delta>0$ small enough, we can define the probability measure 
$$ \nu_\delta = (1+\delta \psi) d\nu,$$
and we have moreover that $L^1(\nu_\delta) = L^1(\nu)$, so that $\phi \in L^1(\nu)$. Since $H(\nu_{\delta} | \mu) = \int \log \big(1 + \delta \psi)e^f\big) d\nu_{\delta}$, we deduce by the Gibbs variational formula \eqref{Gibbs},
$$ H(\nu_\delta | \mu) - H(\nu | \mu) \leq  
 \int \log\big( (1 + \delta \psi)e^f\big) d(\nu_\delta  - \eta). $$
Therefore, dividing \eqref{subdiff} by $\delta$, we get
$$ \int \log (1+\delta \psi) \psi d\nu + \int f \psi d \nu \geq \int \phi \psi d\nu.$$
Taking $\delta \to 0$ we conclude by dominated convergence,
$$ \int f \psi d\nu \geq \int \phi \psi d\nu,$$
for any $\phi$ bounded, measurable such that $\int \psi d\nu =0$. Therefore,
$$ \phi = f + \int (\phi -f)d\nu,$$
$\nu$-almost surely.

Assume now $(i), (ii)$ and that $f$ achieves the equality in \eqref{infconv}. By definition,
$$H(\nu| \mu) = \int f d\nu.$$
As $f$ achieves the equality in \eqref{infconv}, we can write,
$$  \int f d\nu = \int f(x) d\nu(x)  - \int \sup_{x\in \RR^n}\{ f(x) - c(x,y)\} d\tilde \mu(y) \leq \int c(x,y) d\pi(x,y),$$
which proves the second part of the equivalence.
\end{proof}

 As we will see in the sequel, when $\mu$ is a product measure, a transportation-entropy inequality which is saturated by tilts implies by duality a strong integrability inequality for empirical processes.  In this paper we carry out this program in the special case where $\mu$ is the $n$-fold product of a measure supported on $\{-1,1\}$ or of the exponential measure on $\RR_+$.

\subsection{The Gaussian case}\label{Gaussian}
Before going into the investigation of the discrete setting and the case of the exponential measure, we will review what happens in the Gaussian case, which we will regard as a motivational example. We will see how Talagrand's transportation-entropy inequality \cite{Talagrand2} implies a dimension-free mean-field approximation of the free energy of Gibbs measures and a nonlinear large deviation bound.

First, we turn our attention to the mean-field approximation of the free energy. Using Talagrand's transportation-entropy \eqref{TalG} and Proposition \ref{duality}, one gets  that for any measurable function $f : \RR^n \to \RR$,
\begin{equation} \label{infconvG} \log \int e^f d\gamma \leq \int \sup_{h \in \RR^n} \{ f(x+h) - \frac{1}{2}||h||_{\ell^2}^2 \}  d\gamma(x).\end{equation}
Note that replacing $f$ by $\sup_{\lambda \in V} \{\langle \lambda, x\rangle - \frac{1}{2}||\lambda||_{\ell^2}^2\}$, for some countable set $V$, one obtains the strong integrability inequality of Gaussian processes mentioned in \eqref{strongG}. Coming back at the estimation of the free energy, we see that using a change of measure and Jensen's inequality, we have the lower bound
$$ \sup_{h\in \RR^n} \big\{ \int f(x+h) d\gamma(x) - \frac{1}{2}||h||_{\ell^2}^2\big \} \leq \log \int e^f d\gamma.$$
From these two inequalities, we see that the so-called mean-field approximation of the free energy of the Gibbs measure associated with some function $f$ holds as soon as the Gaussian mean-width is small compared to the mean-field approximation. More precisely, we have the following proposition.

\begin{Pro}\label{meanfieldG}
Let $f : \RR^n\to \RR$ be a continuously differentiable function. Then,
\begin{align*}
\sup_{h\in \RR^n} \big\{ \int f(x+h) &d\gamma(x) -\frac{1}{2}||h||_{\ell^2}^2\big\} \leq \log \int e^{f} d\gamma\\
& \leq \sup_{h\in \RR^n} \big\{ \int f(x+h) d\gamma(x) -\frac{1}{2}||h||_{\ell^2}^2\big\} + \sqrt{2}g(V),
\end{align*}
where $V = \nabla f (\RR^n)$ and $g(V)$ is defined in \eqref{Gaussianwidth}.
\end{Pro}
\begin{proof} 
From \eqref{infconvG}, we deduce
\begin{align*}
\log \int e^f d\gamma &\leq  \sup_{h\in \RR^n} \big\{ \int f(x+h) d\gamma(x) -\frac{1}{2}||h||_{\ell^2}^2\big\} \\
&+  \int \sup_{h\in \RR^n}\Big( f(x+h) -\int f(y+h)d\gamma(y)\Big) d\gamma(x).\end{align*}
This last error term can be compared to the Gaussian mean-width of $\nabla f(\RR^n)$, by pulling the integral in $y$ out of the supremum and using the mean-value Theorem, namely,
$$ \int \sup_{h\in \RR^n} \big( f(x+h) - \int f(y+h)d\gamma(y) \big) d\gamma(x)\leq \sqrt{2}g(V), $$
where $V = \nabla f(\RR^n)$. 
\end{proof}

In parallel, Talagrand's transportation-entropy inequality \eqref{TalG} or more strongly the Gaussian isoperimetric inequality can be used to obtain a dimension-free nonlinear large deviations bound. In \cite{LDPWG}, this observation was exploited to derive large deviations principles for a class of functions for which the large deviations are due to translations. In the following proposition, we give a non-asymptotic nonlinear large deviations bound in the Gaussian setting. 

\begin{Pro}\label{NLG}Let $f :  \RR^n \to \RR$ be a continuously differentiable function and denote by $V = \nabla f(\RR^n)$. Define the function,
$$\forall t \in \RR, \ \psi(t) = \inf \big\{ \frac{1}{2}||h||_{\ell^2}^2 : \EE f(X+h)\geq  t, h \in \RR^n\big \}.$$
Let $t \in \RR, \delta>0$. Assume that 
\begin{equation} \label{stricinc} \forall s > t-\delta,\ \psi(s)> \psi(t-\delta).\end{equation}
Let $X$ be a standard Gaussian vector in $\RR^n$. If $g(V) \leq \frac{\delta}{2\sqrt{2}}$, 
then $$\log \PP( f(X) \geq t ) \leq -\psi(t-\delta).$$
\end{Pro}
\begin{Rem}
A close attention to the proof reveals that a weaker sufficient condition is that $\PP(E_{\delta})\geq 1/2$, where
$$ E_{\delta}= \big\{ x : \sup_{||h||_{\ell^2}^2 \leq 2\phi(t-\delta)} | f(x+h) - \EE f(X+h)| <\delta\big\}.$$
This observation can be crucial for certain large deviations problems where the function $f$ does not have a gradient of low complexity in the sense of small Gaussian mean-width but instead in the sense that $E_\delta$ is a typical set. An example of such a large deviation problem is given by the traces of power of Gaussian Wigner matrices which was studied in \cite{LDPtr}.

\end{Rem}
\begin{proof}Let $Y$ be a standard Gaussian random variable independent of $Y$. As a consequence of the assumption \eqref{stricinc}, we claim that
\begin{equation} \label{bound} \PP\big( f(X) \geq t \big) \leq \PP\Big( \inf_{ || h||_{\ell^2}^2 \leq 2 \psi(t-\delta)}\big( f(X) - \EE f(Y+h) \big)\geq\delta \Big).\end{equation}
Indeed, if $||h||_{\ell^2}^2 \leq 2 \psi(t-\delta)$, then by the definition of $\psi$, 
$$ \psi\big( \EE f(Y+h) \big) \leq \psi(t-\delta).$$
Therefore $\EE f(Y+h) \leq t-\delta$ by \eqref{stricinc}. Define the set
$$ E_{\delta}= \big\{ x : \sup_{||h||_{\ell^2}^2 \leq 2\psi(t-\delta)} (f(x+h) - \EE f(Y+h))<\delta\big\}.$$
With this notation, one can observe that,
$$ \PP\Big( \inf_{ || h||_{\ell^2}^2 \leq 2 \psi(t-\delta)} \big( f(X) - \EE f(Y+h)\big) \geq\delta \Big) \leq \PP\big( X \notin E_{\delta} + \sqrt{2\psi(t-\delta)} B_{\ell^2}\big).$$
The Gaussian isoperimetric inequality entails that the Gaussian measure has normal concentration by \cite[Corollary 2.6]{Ledouxmono}, which means that,
$$\PP\big( X \notin E_{\delta} + \sqrt{2\psi(t-\delta)} B_{\ell^2}\big) \leq e^{-\psi(t-\delta)},$$
as soon as $\PP(X \in  E_\delta) \geq 1/2$. The mean-value Theorem and Markov's inequality yield,
$$ \PP(X \notin  E_\delta) \leq \frac{1}{\delta} \EE   \sup_{x \in \RR^n} \langle \nabla f(x), X-Y\rangle,$$
which concludes the proof.
\end{proof}

Let us make some closing remarks about the specificity of the Gaussian setup. 
For a convex function $\Lambda : \RR^n \to \RR \cup\{+\infty\}$, one defines its \textit{Legendre transform} by the formula,
$$ \forall x \in \RR^n, \ \Lambda^*(x) = \sup_{\xi \in \RR^n} \big\{ \langle \xi,x\rangle - \Lambda(\xi) \big\}.$$
We will also denote by $\Gamma(\RR^n)$ the set of convex functions $\RR^n \to \RR\cup\{+\infty\}$ which are lower semi-continuous, and are proper, namely their domain is nonempty.
 In general, one can show using a small modification of \cite[Remark 2.12]{LW} the following fact on cost functions of transportation-entropy inequalities.

\begin{Fact} Let $\mu$ be a probability measure on $\RR^n$ with mean $0$ which  satisfies a transportation-entropy inequality \eqref{transpoin} with cost function $c$ of the form,
$$ c : (x,y) \mapsto \alpha(x-y),$$
for some function $\alpha \in \Gamma(\RR^n)$. Then,
$$\forall x \in \RR^n, \  \alpha(x) \leq \Lambda^*_\mu(x).$$ 
\end{Fact}
\begin{proof}
By Proposition \ref{duality}, we know that for any measurable function $f : \RR^n \to \RR$, such that $f^c \in L^1(\mu)$,
$$ \log \int e^{f}d\mu \leq \int \sup_{x \in \RR^n} \{ f(x)-\alpha(x-y) \}  d\mu(y).$$ 
Testing the above inequality for linear forms we get
$$\forall \theta \in \RR^n, \  \Lambda_\mu(\theta) \leq \alpha^*(\theta). $$
As $\alpha \in \Gamma(\RR^n)$, we can conclude using \cite[ Theorem 4.21]{Clarke} that $\alpha \leq \Lambda^*_\mu$.
 
\end{proof}

Therefore, among the cost functions of the form $c(x,y) =\alpha(x-y)$, with $\alpha \in \Gamma(\RR^n)$, the best cost function one can expect is
$$ (x,y) \mapsto \Lambda^*_\mu(x-y).$$
 Note that Talagrand's result \eqref{TalG} gives exactly that the Gaussian measure satisfies a transportation-entropy inequality with the above cost function.
As we now show,  the Gaussian measure is \textit{the only} measure with this property.

\begin{Fact}
Let $\mu$ be a probability measure on $\RR^n$ with mean $0$ such that the domain of $\Lambda_\mu$ has nonempty interior and $\Lambda_\mu^*$ is strictly convex. If $\mu$ satisfies the transportation-entropy inequality with cost function
$$ (x,y) \mapsto \Lambda_\mu^*(x-y),$$
then $\mu$ is a Gaussian measure.
\end{Fact}

\begin{Rem}
By \cite[Theorem 4.1]{Rockafeller}, a sufficient condition for $\Lambda^*_\mu$ to be strictly convex is that the support of $\mu$ is not included in a hyperplane and $\Lambda_\mu$ is essentially convex, that is,  its domain, denoted by $\mathcal{D}_\Lambda$, has nonempty interior, $\Lambda_\mu$ is differentiable on the interior of its domain $\mathcal{D}_\Lambda^\circ$, and steep, that is, for any $\xi_k \in \mathcal{D}_\Lambda^\circ$ such that $\xi_k \to \xi \in \partial \mathcal{D}_\Lambda$ when $k \to +\infty$, we have
$$ || \nabla \Lambda (\xi_k) ||_{\ell^2} \underset{k\to+\infty}{ \longrightarrow} +\infty.$$
\end{Rem}

\begin{proof}  
We write as a short-hand $\Lambda$ instead of $\Lambda_\mu$ and $\mathcal{D}_\Lambda$ its domain. Let $\theta \in \mathcal{D}_{\Lambda}^{\circ}$, and define the probability measure
$$\mu_{\theta} = e^{\langle x, \theta \rangle - \Lambda(\theta)} d\mu(x).$$
As $\Lambda^*$ is a convex function, it is continuous on the interior of its domain. Therefore, the function
$$  c: (x,y) \mapsto \Lambda^*(x-y),$$
is lower semi-continuous. By \cite[Theorem 4.1]{Villani}, we know that there exists a coupling $\pi$ between $\mu_{\theta}$ and $\mu$, such that
$$  \mathcal{W}_c(\mu_\theta, \mu) = \int \Lambda^*(x-y) d\pi(x,y).$$
Assume $\mu$ satisfies the transportation-entropy with cost function $c$. As 
$$ H(\mu_\theta | \mu) = \Lambda^*(\nabla \Lambda(\theta)),$$
we have 
$$ \int \Lambda^*(y-x) d\pi(x,y) \leq \Lambda^*(\nabla \Lambda(\theta)).$$
But by convexity of $\Lambda^*$, we get,
$$ \Lambda^*(\nabla \Lambda(\theta)) = \Lambda^*\big( \int y d\mu_{\theta}(y) \big)   \leq \int \Lambda^*(y-x) d\pi(x,y) \leq \Lambda^*(\nabla \Lambda(\theta)).$$
As $\Lambda^*$ is strictly convex, the equality in Jensen's inequality yields that $\mu_{\theta}$ is the push forward of $\mu$ by a translation. Since the mean of $\mu_\theta$ is $\nabla \Lambda(\theta)$ and the one of $\mu$ is $0$, $\mu_\theta$ is the push-forward of $\mu$ by the map $x \mapsto x + \nabla \Lambda(\theta)$. Comparing the log-Laplace transforms of $\mu_\theta$ on one hand, and the one of $\mu$ pushed-forward by $x \mapsto x +\nabla \Lambda(\theta)$, we find that $\Lambda$ satisfies the following functional equation:
\begin{equation} \label{eqfunc} \forall \theta \in \mathcal{D}_\Lambda, \xi \in \RR^n, \ \Lambda(\xi+\theta) =  \Lambda(\theta) + \Lambda(\xi) + \langle \xi,  \nabla \Lambda(\theta) \rangle.\end{equation}
From this equation, we see that if $\xi, \theta \in \mathcal{D}_\Lambda$, then $\xi+ \theta \in \mathcal{D}_\Lambda$. As the interior of $\mathcal{D}_\Lambda$ is nonempty, we must have $ \mathcal{D}_\Lambda = \RR^n$. 
Differentiating \eqref{eqfunc} with respect to $\xi$, we get
$$ \forall \theta, \xi \in \RR^n, \ \nabla \Lambda(\xi+\theta) = \nabla \Lambda(\xi) + \nabla \Lambda(\theta).$$
As $\nabla \Lambda$ is continuous, the above equation implies that $\nabla \Lambda$ is a linear function. Thus, $\Lambda$ is a quadratic form and $\mu$ is a Gaussian measure.
\end{proof}

Even though the Gaussian measure is the only measure $\mu$ to satisfy a transportation-entropy inequality with cost function $\Lambda^*_\mu(x-y)$, 
it has been shown in \cite{LW} that up to some universal constant $\beta>0$, any symmetric log-concave product measure $\mu$ on $\RR$ satisfies a transportation-entropy inequality with cost function $\Lambda^*_\mu\big( \frac{x-y}{\beta}\big)$.  A similar result has been proven for measures with log-concave tails in \cite{SST} and weak transport-entropy inequalities. However, following the argument of the proof of Proposition \ref{meanfieldG}, we see that it entails that for any continuously differentiable $f : \RR^n \to \RR$, and $\mu$ a symmetric log-concave product measure, 
\begin{align*}
 \log \int e^f d\mu & \leq \sup_{h \in \RR^n}\big\{ \int f(x+h) d\mu(x)  - \Lambda^*_{\mu}\Big(\frac{h}{\beta}\Big) \big\}\\&  + \int \sup_{z \in \RR^n} \langle \nabla f(z), x-y \rangle d\mu(x) d\mu(y),\end{align*}
 Thus, it yields a multiplicative error from the true entropic cost one expects. But in the applications we have in mind, it will be important for us to have the best constant, that is $\beta = 1$, in order to produce a matching upper bound, so that we cannot rely on the mentioned results.

\subsection{The exponential measure}

The moral we deduce from the Gaussian case is that we have to look for cost functions beyond the ones of the form $(x,y)\mapsto \alpha(x-y)$, in order to hope for transportation-entropy inequalities to be saturated by tilted measures. 
In the case of the exponential measure, we consider the cost function,
\begin{equation} \label{defcoutexp1}\forall x,y \in \RR^n, \ c(x,y) = \sum_{i=1}^n  y_i \Lambda^*_{\eta}\Big( \frac{x_i}{y_i}\Big).\end{equation}
The form of this cost function can be explained by the natural coupling of all the tilts $(\eta_\lambda)_\lambda$ of the exponential measure, where 
$$  \forall \lambda >0, \ \eta_{\lambda} = \Car_{x\geq 0} \lambda e^{-x/\lambda} dx.$$
There is a simple way to transport $\nu$ onto $\nu_{\lambda}$ by the map $x \mapsto \lambda x$. This fact explains the shape of the cost function \eqref{defcoutexp1} as essentially a function of the ratio $y/x$. 
We now give a proof of Proposition \ref{transpoexp}.
\begin{proof}[Proof of Proposition \ref{transpoexp}] 
By the tensorization property of the transportation-entropy inequalities (see Proposition \ref{tensor}), it is sufficient to prove the statement for $n=1$. Let $\nu$ be a probability measure on $\RR_+$. 
Let $\tilde \nu$ and $\tilde \eta$ be the push-forward of respectively $\nu$ and $\eta$ by the map $x \mapsto \log x$. Note that, 
$$ \tilde \eta = e^{-\xi(x)} dx,$$
with $\xi(x) = e^x-x$, which is a convex function. 
From \cite{Talagrand2} we know that $\tilde \eta$ satisfies a transportation-entropy inequality with cost function $\tilde{c}$ defined by,
$$ \forall x,y \in \RR,\ \tilde{c}(x,y) =  \xi(x) - \xi(y) -\xi'(y)(x-y),$$
that is,
$$ \tilde c(x,y) = e^y \Lambda^*_\nu(e^{x-y}).$$ 
But on one hand,
$$ H(\tilde \nu | \tilde \eta) = H( \nu | \eta).$$
On the other hand,
$$ \mathcal{W}_{\tilde c}(\tilde \nu , \tilde \eta) = \mathcal{W}_c(\nu, \eta),$$
which gives the first claim.

It only remains to prove that if $\nu$ is a tilt of $\eta^n$ then it achieves the equality in the transportation-entropy inequality of Proposition \ref{transpoexp}. Let $\lambda = (\lambda_1,\ldots,\lambda_n)$ with $\lambda_i >0$. Denote by $\eta_\lambda = \eta_{\lambda_i} \otimes \ldots \otimes \eta_{\lambda_n}$. Let $\pi$ be a coupling between $\eta_\lambda$ and $\eta^n$. By Jensen's inequality we have,
$$ \sum_{i=1}^n \Lambda^*_\eta\big( \int x_i d\eta_{\lambda_i}(x_i)\big) \leq \int \sum_{i=1}^n y_i \Lambda^*_\eta\Big( \frac{x_i}{y_i}\Big) d\pi(x,y).$$
On the other hand,
$$ H( \eta_{\lambda}| \eta^n) = \sum_{i=1}^n \Lambda^*_\eta(\lambda_i).$$
Thus, 
$$ \mathcal{W}_c(\eta_\lambda| \eta^n) \geq H(\eta_\lambda| \eta^n),$$
which ends the proof of the equality case.
\end{proof}

Using duality, we can now give a proof of Proposition \ref{intexpo} on the integrability of empirical exponential processes. 

\begin{proof}[Proof of Proposition \ref{intexpo}]Let $V \subset \RR^n$ be a countable subset. Define the function $g$ by,
$$\forall x \in \RR^n, \  g(x) = \sup_{\xi \in V } \{ \langle \xi, x\rangle -\Lambda_{\eta^n}(\xi)\}.$$
By Propositions \ref{transpoexp} and \ref{duality} we have,
$$ \log \int e^{g}d\nu \leq \int \sup_{y \in \RR_+^n} \{ g(y) - c(x,y) \}d\mu(x),$$
where $c$ is defined in \eqref{defcoutexp1}.
But,
$$ \sup_{y\in \RR_+^n} \{ g(y) - c(x,y) \} = \sup_{t \in \RR_+^n} \{ g(tx)-\langle   x, \Lambda^*_\eta(t)\rangle \},$$
where $tx = (t_1 x_1,\ldots ,t_nx_n)$ and $\Lambda_\eta^*(t) = (\Lambda^*_\eta(t_1),\ldots,\Lambda_\eta(t_n))$.
Therefore,
$$ \sup_{y\in \RR_+^n} \{ g(y) -c(x,y) \} = \sup_{\xi \in V} \sup_{t \in \RR_+^n} \{ \langle \xi, tx\rangle-\langle x,\Lambda^*_{\eta}(t)\rangle -\Lambda_{\eta^n}(\xi) \}.$$
Fix $\xi \in V$. We have
\begin{align*} \sup_{t \in \RR_+^n} \{ \langle \xi, tx\rangle-\langle x, \Lambda^*_\eta(t)\rangle \}& = \sum_{i=1}^n \sup_{t>0} (t \xi_i-\Lambda^*_\eta(t)) x_i \\
& = \sum_{i=1}^n  \Lambda_{\eta}(\xi_i) x_i,
\end{align*}
where we used the fact that $\Lambda_{\eta}$ is the Legendre transform of $\Lambda^*_{\eta}$.
\end{proof}

\subsection{The discrete hypercube}
Let $p\in(0,1)$ and $\mu_p= (1-p) \delta_{-1} + p \delta_{1}$. Note that the tilts of $\mu_p^n$ are exactly the product probability measures on  $\{-1,1\}^n$. One of the difficulties in finding a transportation-entropy inequality on the discrete hypercube which is saturated by product measures comes from the fact that the measure $\mu_p^n$ does not carry enough information in order to sample from it all product measures on $\{-1,1\}^n$. This assessment brings us to the conclusion that one has to enrich the background measure, and consider not a transportation cost between a given probability measure $\nu$ on $\{-1,1\}^n$ and $\mu_p^n$, but a transportation cost between $\nu$ and the uniform measure on $[-1,1]^n$. In this section, we prove the transportation-entropy inequality of Proposition \ref{transpo}.

By the tensorization property of transportation-entropy inequalities, which we recalled in Proposition \ref{tensor}, we only have to prove Proposition \ref{transpo} for $n=1$. This is the content of the following lemma.

\begin{Lem}\label{init}
 For any $\nu$ probability measure on $\{-1,1\}$,
\begin{equation} \label{transpo1}\mathcal{W}_{w_p}(\nu, \mathcal{U}) = H(\nu | \mu_p),\end{equation}
where $w_p$ is defined in \eqref{defw}, and $\mathcal{U}$ denotes the uniform measure on $[-1,1]$.
\end{Lem}
\begin{proof}
Let $h$ denote the mean of $\nu$. Let $\pi$ be the law of
$$ (U,\sg(h-U)),$$
where $U$ is uniformly distributed on $[-1,1]$, and $\sg(x) = 1$ if $x\geq 0$ and $-1$ otherwise.
By definition of $\pi$,
$$ \int w_p(x,u) d\pi(x,u) = 2\EE  |I'_p(U)| \Car_{U\in (h,h_0)} = I_p(h),$$
using the fact that $I_p(h_0) = 0$.
Besides,
$$ H(\nu| \mu_p) = I_p(h),$$
which proves the inequality
$$ \mathcal{W}_{w_p}(\nu,\mathcal{U})\leq H(\nu|\mu_p).$$
To prove the equality, we will prove that equality is achieved in the inf-convolution inequality with cost function $w_p$, where  we set $w_p = +\infty$ on $\RR\setminus \{-1,1\}$, and use Proposition \ref{duality}.
Let $t \in \RR$, and define 
$$ Y_t =  \max_{x \in \{-1,1\}} \{tx -w_p(x,U) \},$$
with $U$ uniformly sampled in $[-1,1]$. We need to prove that 
\begin{equation} \label{optimexp}\EE Y_t = \Lambda_{p}(t),\end{equation}
where we use $\Lambda_p$ as a short-hand for $\Lambda_{\mu_p}$. We have,
\begin{align*}
Y_t & = \max\big( t -2I_p'(U)\Car_{U>h_0}, -t +2I_p'(U)\Car_{U <h_0}\big)\\
& = \max\big( t-2I_p'(U), -t \big)\Car
_{U>h_0} + \max\big( t, -t+2I_p'(U)\big)\Car_{U<h_0}.
\end{align*}
There are two cases to consider. First we assume that $h_0 \leq \Lambda_p'(t)$. Observe that $I_p$ is the Legendre transform of $\Lambda_p$. Therefore, $I_p'$ and $\Lambda_p'$ are inverse functions. 
We can write,
$$Y_t  = \big(t-2I_p'(U)\big)\Car_{h_0 \leq U \leq \Lambda_p'(t)} -t \Car_{U>\Lambda_p'(t)}+ t \Car_{U<h_0}.$$
Thus, 
$$ \EE Y_t = \frac{t}{2} ( \Lambda_p'(t) - h_0) - I_p(\Lambda'(t))  - \frac{t}{2}(1- \Lambda_p'(t))+\frac{t}{2}(h_0+1),$$
where we used the fact that $I_p(h_0) = 0$. Therefore, 
$$ \EE Y_t = t \Lambda_p'(t) -I_p(\Lambda_p'(t)) = \Lambda_p(t), $$
since $I_p$ is the Legendre transform of $\Lambda_p$. If $h_0 >\Lambda_p'(t)$, we get 
$$Y_t  = -t \Car_{U\geq h_0} + t \Car_{U \leq \Lambda_p'(t)} +\big(-t+2 I_p'(U)\big)\Car_{\Lambda_p'(t) < U <h_0},$$ 
which yields similarly $\EE Y_t = \Lambda_p(t)$.

\end{proof}

\section{Strong integrability of  empirical processes}
In this section, we show how transportation-entropy inequalities which are saturated by tilts implies a strong integrability inequality for empirical processes. 

\begin{Pro}\label{strongintgene}
Let $\mu$ be a probability measure on $\RR$ with support included in $[-1,1]$. 
Let $\tilde \mu$ be a probability measure on $\RR$ and let $w : \RR\times \RR \to [0,+\infty]$ be a lower semi-continuous function such that 
$$ \forall \nu \in \mathcal{P}(\RR), \ \mathcal{W}_w(\nu,\tilde \mu) \leq H(\nu | \mu),$$
and equality holds for the tilts of $\mu$. 
Then, for any countable subset $V \subset \RR^n$,
$$ \log \int e^{\sup_{\xi \in V} \{ \langle x,\xi \rangle - \Lambda_{\mu^n}(\xi)\} } d\mu^n(x) \leq \kappa b(V),$$
where $\kappa$ is a universal constant, and $b(V)$ is defined in \eqref{Rademwidth}.
\end{Pro}
Combining Proposition \ref{strongintgene} with Proposition \ref{transpo}, we obtain the result of Proposition \ref{strongintB}.

\begin{proof} As $\mu$ has its support included in $[-1,1]$, we can assume without loss of generality that $w(x,y) = +\infty$ whenever $x \notin [-1,1]$. By the tensorization property of transport-entropy inequality (see Proposition \ref{tensor}), $\mu^n$ satisfies the transportation-entropy inequality,
$$ \forall \nu \in \mathcal{P}(\RR^n), \ \mathcal{W}_w(\nu,\tilde \mu^n) \leq H(\nu|\mu^n),$$
where
$$ \forall x,y \in \RR^n, \ w(x,y ) := \sum_{i=1}^n w(x_i,y_i).$$
 Moreover, from the Fact \ref{tensoreq}, equality holds for the tilts of $\mu^n$.

Let $V$ be a countable subset of $\RR^n$.  By Proposition \ref{duality}, we have
$$ \int e^{\sup_{\xi \in V}\{ \langle \xi, x\rangle - \Lambda_{\mu^n}(\xi) \}} d\mu^n(x) \leq \int \phi d\tilde \mu^n,$$
where 
$$ \forall y \in \RR^n, \ \phi(y) = \sup_{\xi \in V}\sup_{ y\in [-r,r]^n}  \{ \langle \xi, x\rangle -w(x,y) -\Lambda_{\mu^n}(\xi)\}.$$
Thus, it remains to compute the expectation of a supremum of a certain empirical process. To this end, we will use the characterization of the boundedness of Bernoulli processes, proven by Bednorz and Lata\l{}a \cite{BL}. We start by showing that this process has sub-Gaussian increments with variance factor given by the $\ell^2$-norm. We write in probabilistic notation, 
$$ \int \phi d\tilde \mu^n = \EE \sup_{\xi\in V} Z_{\xi},$$
where for any $\xi \in \RR^n$,
$$ Z_{\xi} = \sum_{i=1}^n \big(T_{\xi_i}- \Lambda_\mu(\xi_i)\big),$$
with $T_{\xi_i} =  \max_{x \in [-r,r]} \{ \xi_i x -w(x,Y_i) \}$, and $Y=(Y_1,\ldots ,Y_n)$ is sampled according to $\tilde \nu^n$. 
%

Let $\xi,\zeta \in \RR^n$.  For any $i\in\{1,\ldots,n\}$, we have
$$  T_{\xi_i}-T_{\zeta_i} \leq  |\xi_i- \zeta_i|.$$
The fundamental fact about the process $Z_\xi$ is that it is centered. This is due to the fact that equality in the transportation-entropy inequality with cost function $w$ holds for tilts of $\mu$. Indeed, by Proposition \ref{duality}, we deduce that equality holds in the corresponding inf-convolution inequality for linear forms, which exactly says that,
$$ \EE T_{\xi_i} = \Lambda_\mu(\xi_i).$$
Therefore, we can write
 $$Z_{\xi} - Z_{\zeta} = \sum_{i=1}^n (\Delta_i - \EE\Delta_i),$$
where $\Delta_i$ are independent, and $| \Delta_i| \leq |\xi_i - \zeta_i|$. Thus, by Hoeffding's inequality (see \cite[Theorem 2.8]{BLM}), for any $t>0$,
\begin{equation} \label{incre} \PP\big( |Z_\xi - Z_\zeta| >t \big) \leq \exp\Big( - \frac{t^2}{2 ||\xi - \zeta ||_{\ell^2}^2}\Big).\end{equation}
Therefore, if $V\subset V_1 +V_2$, then for any $\xi \in V$, $\xi = \xi_{1} + \xi_2$, with $\xi_i \in V_i$, 
\begin{align*}
 Z_{\xi}& = \sup_{x \in [-r,r]^n}\{ \langle\xi,x \rangle -w(x,Y) \} -\Lambda_{\mu^n}(\xi)\\
& \leq  \sup_{x \in [-r,r]^n}\{ \langle\xi_2,x \rangle -w(x,Y) \} -\Lambda_{\mu^n}(\xi_2)\\
& + ||\xi_1||_{\ell^1} +  \Lambda_{\mu^n}(\xi_2)-\Lambda_{\mu^n}(\xi_1+\xi_2).
\end{align*}
As the support of $\mu$ is included in $[-1,1]$, $\Lambda_{\mu^n}$ is $1$-Lipschitz with respect to the $\ell^1$-norm. Therefore, 
$$  Z_{\xi} \leq Z_{\xi_2} + 2 ||\xi_1||_{\ell^1}.$$
But, from the incremental property \eqref{incre} and the Majorization Theorem (see \cite[Theorem 12.16]{LT}), we get 
$$ \EE \sup_{\xi \in V_2 } Z_{\xi} \leq L g(V_2),$$
where $L$ is a numerical constant, and $g(V_2)$ is the Gaussian width of $V_2$, that is,
$$ g(V_2) = \EE \sup_{\xi \in V_2} \langle \xi, \Gamma \rangle,$$
where $\Gamma$ is a standard Gaussian random variable in $\RR^n$. Therefore,
$$ \EE \sup_{\xi \in V} Z_{\xi} \leq r \inf\big\{ 2 \sup_{\xi \in V_1} ||\xi||_{\ell^1} + L g(V_2) : V \subset V_1 +V_2 \big\}.$$
We know from the characterization of the boundedness of Bernoulli processes \cite{BL} that there exists a numerical constant $C>0$ such that,
$$\inf\big\{ \sup_{\xi \in V_1} ||\xi||_{\ell^1} + g(V_2) : V \subset V_1 +V_2 \big\} \leq C  b(V),$$
which gives the claim.

\end{proof}
Using the same comparison arguments between supremum of empirical processes, we obtain the following lemma which enables to estimate the difference between extensions $\{-1,1\}^n$ to $[-1,1]^n$ of functions with a low complexity set of gradient. 
\begin{Lem} \label{exten}
Let $f : \RR^n \to \RR$ be a continuously differentiable function and $g$ the harmonic extension of $f_{|\{-1,1\}^n}$ to $[-1,1]^n$ defined as in remark \ref{harmoext}. Then,
$$ \sup_{ y \in [-1,1]^n } \big( f(y) -g(y) \big) \leq \kappa b(V),$$
where $\kappa$ is a numerical constant, $V = \nabla f( [-1,1]^n)$ and $b(V)$ is defined as in \eqref{Rademwidth}.
\end{Lem}
\begin{proof}
Let $y \in [-1,1]^n$. By definition,
$$  f(y)  - g(y)=  \EE f(y) - f(X_y),$$
where $X_y$ is a  random vector in $\{-1,1\}^n$ with independent coordinates and mean $y$. 
By the mean-value Theorem, we have
$$ f(y) - g(y)\leq \EE \sup_{ \xi \in V} \langle \xi, X_y-y \rangle.$$
Repeating the argument of the proof of Proposition \ref{strongintgene}, we see that the characterization of boundedness of Bernoulli processes \cite{BL} entails that there exists a universal constant $\lambda >0$, such that
$$ f(y)- g(y) \leq \lambda b(V).$$

\end{proof}

\section{The mean-field approximation}
Building on the strong integrability of Bernoulli empirical processes, we give here a proof of the mean-field approximation of the free energy of Gibbs measures on the discrete hypercube stated in Theorem \ref{meanfield}.

In a first step we identify the error term induced by the mean-field approximation with the help of the following lemma. A proof of this result can be found in \cite[Proposition 1.1]{NL}.

\begin{Lem}\label{meanfieldstep}
Let $\mu$ be a compactly supported probability measure on $\RR^n$. Denote by $K$ the convex hull of its support. For any $f :\RR^n \to \RR$ continuously differentiable,
$$ \log \int e^f d\mu \leq \sup\{ f- \Lambda^*_\mu\}+\log \int e^{\sup_{\xi \in V} \{ \langle\xi, x\rangle - \Lambda_\mu(\xi) \}} d\mu(x),$$
where $ V$ is the convex hull of $\nabla f(K)$, and $\Lambda^*_\mu$ denotes the Legendre transform of $\Lambda_\mu$.
\end{Lem}

Combining Lemma \ref{meanfieldstep} with Propositions \ref{transpo} and \ref{strongintgene}, we obtain Theorem \ref{meanfield}. In fact, we have the following more general result which states that a dimension-free mean-field approximation holds as soon as a transportation-entropy inequality is saturated by tilts exists.

\begin{Pro}
Let $\mu$ be a probability measure on $\RR$ with support included in $[-1,1]$. 
Assume there exist $\tilde \mu$ a probability measure on $\RR$ and $w : \RR\times \RR \to [0,+\infty]$ a lower semi-continuous function such that 
$$ \forall \nu \in \mathcal{P}(\RR), \ \mathcal{W}_w(\nu,\tilde \mu) \leq H(\nu | \mu),$$
and equality holds for the tilts of $\mu$. Then, there exists a universal constant $\kappa>0$ such that for any $f : \RR^n \to \RR$ continuously differentiable,
$$ \log \int e^f d\mu \leq \sup\{ f- \Lambda^*_{\mu}\} + \kappa b(V),$$
where $V = \nabla f([-1,1]^n)$ and $b(V)$ is defined in \eqref{Rademwidth}.
\end{Pro}

\section{Proof of Theorem \ref{NL}}
Contrary to the nonlinear large deviations bounds shown in the previous works \cite{CD}, \cite{Yan}, and \cite{Eldan}, the proof of Theorem \ref{NL} will not rely on the computation of exponential moments of functions with a low complexity set of gradients. Instead, we will show as a first step that one can  reformulate the deviations of $f(X)$ in terms of the deviations of the process $(\langle \theta \xi, X\rangle - \Lambda_\mu(\theta \xi))_{ \xi \in V, \theta>0}$.  Then, we will use the strong integrability inequality of Bernoulli processes of Proposition \ref{strongintB} to control the deviations of the latter process.

In the next lemma, we relate the deviations of $f(X)$ and of the process $(\langle \theta \xi, X\rangle - \Lambda_\mu(\theta \xi))_{ \xi \in V, \theta>0}$. We state it in the general setting where $X$ is distributed to a compactly supported measure.
\begin{Lem}
Let $\mu$ be a compactly supported probability measure on $\RR^n$, whose support is not included in a hyperplane. Denote by $K$ the convex hull of its support. Let $f :\RR^n \to \RR$ be a continuously differentiable function and let $W$ be the convex hull of $\nabla f(K)$. Define the function,
$$ \forall t \in \RR, \ \phi(y) = \inf \big\{ \Lambda^*_\mu(y) : f(y) \geq t \big\}.$$ 
 Let $t\in \RR$ and $\delta>0$. Assume that 
$$ \forall s > t-\delta, \ \phi(t-\delta) < \phi(s).$$
For any $x \in K$,
$$  f(x) \geq t \Longrightarrow  \sup_{\xi \in W \atop 0\leq \theta \leq \theta_0} \{ \langle \theta \xi, x\rangle - \Lambda_\mu(\theta \xi) - \theta \delta\} \geq \phi(t-\delta),$$
where $\theta_0 = \Lambda^*_\mu(x) / \delta$.
\end{Lem}
\begin{proof}
Let $x \in K$ such that and $f(x) \geq t$. Arguing as in the proof of Proposition \ref{NLG}, we have
$$ \inf_{\Lambda^*_\mu(y) \leq \phi(t-\delta)} \big(f(x) - f(y)\big) \geq \delta.$$
By the mean value Theorem, we deduce that
$$ \inf_{\Lambda^*_\mu(y) \leq \phi(t-\delta)} \sup_{\xi \in W} \langle \xi, x-y \rangle \geq \delta,$$
which means,
$$ \inf_{h_W(z) \leq \delta} \Lambda^*_\mu(x-z) \geq \phi(t-\delta),$$
where $h_W$ denotes the support function of $W$, namely, $h_W(z) = \sup_{\xi \in W} \langle \xi, z\rangle$. Note that $\Lambda^*_\mu = +\infty$ on $K^c$ by the Hahn-Banach Theorem. Moreover, $\Lambda_\mu^*$ is lower semi-continuous as it is a Legendre transform. Thus $\Lambda^*_\mu$ has compact level sets. As $\{ h_W \leq \delta\}$ is closed, we deduce that the infimum of $\Lambda^*_\mu(x-.)$ on $\{ h_W \leq \delta\}$ is achieved at some $z_*$.

Since $\Lambda^*_\mu$ and $h_W$ are both convex functions, we deduce by Kuhn-Tucker Theorem (see \cite[Theorem 9.4]{Clarke}) that there exists $(\eta, \theta)\neq (0,0)$ with $\eta \in\{0,1\}$ and $\theta \geq 0$, such that $\theta(h_W(z_*) -\delta) = 0$, and
\begin{equation}  \label{inf} \eta \Lambda^*_\mu(x-z_*) = \inf \big\{ \eta \Lambda^*_\mu(x-z) + \theta (h_W(z)-\delta) : z \in \RR^n \big\}.\end{equation}
Evaluating the function on the right-hand side at $z=0$, we see that the non-triviality condition $(\eta,\theta) \neq (0,0)$ implies that $\eta = 1$. Moreover, 
$$ \Lambda^*_\mu(x) - \theta \delta \geq \Lambda^*_\mu(x-z_*) \geq 0.$$
Thus, $\theta \leq \Lambda^*_\mu(x)/\delta$. 
 As $\Lambda^*_\mu = +\infty$ on $K^c$, the infimum in \eqref{inf} can be restricted to $x-K$. Using the Minimax Theorem (see \cite[Theorem 4.36]{Clarke}), we obtain,
$$ \inf \big\{ \Lambda^*_\mu(x-z) + \theta h_W(z) : z \in x-K \big\}= \sup_{\xi \in W} \inf_{z \in \RR^n} \big\{ \Lambda^*_\mu(x-z) +\theta \langle \xi, z\rangle \big\}.$$
We can identify this later infimum using the fact that $\Lambda_\mu$ is the Legendre transform of $\Lambda^*_\mu$ by \cite[Theorem 4.21]{Clarke},
$$\inf_{z \in \RR^n} \big\{ \Lambda^*_\mu(x-z) +\theta \langle \xi, z\rangle \big\} = \langle \theta \xi, x\rangle - \Lambda_\mu(\theta \xi),$$
which ends the claim.
 \end{proof}

We now come back to the proof of Theorem \ref{NL}.   We have
$$ \forall  x \in [-1,1]^n, \ I_p(x) \leq n \log\Big(\frac{1}{p(1-p)}\Big).$$
Denoting by $\Lambda_p$ the log-Laplace transform of $\mu_p^n$, and using the preceding lemma, we get,
$$\PP(f(X) \geq t) \leq \PP\Big( \sup_{\xi \in W \atop 0\leq \theta \leq \theta_0} \{ \langle \theta \xi, X\rangle - \Lambda_{p}(\theta \xi) - \theta \delta\} \geq \phi_p(t-\delta)\Big),$$
with $\theta_0 = -n \log(p(1-p))/\delta$.
We now perform a net argument on $\theta$. Let $\mathcal{D}$ be a $1/(2\sqrt{n}L)$-net of the interval $[0,\theta_0]$, where 
$$ L = \sup_{ x\in K} ||\nabla f (x) ||_{\ell^2} = \sup_{ \lambda \in W} ||\lambda ||_{\ell^2},$$
One can find a net $\mathcal{D}$ such that,
$$ |\mathcal{D}| \leq \frac{4 n\sqrt{n} \kappa |\log(p(1-p))|L }{\delta}.$$
 For $X\in \{-1,1\}^n$ fixed, define the function 
$$G : \xi \in \RR_+  \mapsto \sup_{\xi \in W} \{ \langle \theta\xi, X\rangle - \Lambda_p(\theta \xi) - \theta \delta\}.$$
We claim that for any $\theta' \leq \theta$,
\begin{equation} \label{liponeside} G(\theta) - G(\theta') \leq 2(\theta -\theta')L\sqrt{n}.\end{equation}
Indeed, there is some $\xi \in V$ such that,
$$ G(\theta) -G(\theta') \leq (\theta-\theta')\langle \xi, X\rangle   -\Lambda_p(\theta \xi) + \Lambda_p(\theta'\xi) - (\theta-\theta')\delta.$$
By convexity,
$$ G(\theta) -G(\theta') \leq (\theta-\theta')\langle \xi, X - \nabla \Lambda_p(\theta' \xi) \rangle.$$
As $|| \xi||_{\ell^2} \leq L$ and $\nabla \Lambda_p(\theta' \xi) \in [-1,1]^n$, we get the claim \eqref{liponeside}.

Thus, using a union bound, we get,
\begin{align*}
 \PP& \Big( \sup_{ \xi \in  W \atop 0 \leq \theta \leq \theta_0} \{ \langle \theta \xi, X\rangle - \Lambda_p(\theta\xi) - \theta \delta\} \geq \phi_p(t-\delta)\Big)   \\
&\leq \sum_{\theta \in \mathcal{D}} \PP\Big( \sup_{ \xi \in  W}  \{ \langle \theta \xi , X\rangle - \Lambda_p(\theta\xi) - \theta \delta\} \geq \phi_p(t-\delta)-1\Big).\end{align*}
Now, fix $\theta \in \mathcal{D}$. By Chernof's inequality, we have
\begin{align*}
\log \PP\Big( \sup_{ \xi \in  W}  \{ \langle \theta \xi , X\rangle - \Lambda_p(\theta\xi) - \theta \delta\}& \geq \phi_p(t-\delta)-1\Big) \leq  - \phi_p(t-\delta) +1 \\
&+ \log \EE  e^{\sup_{ \xi \in  W}  \{ \langle \theta \xi , X\rangle - \Lambda_p(\theta\xi) \}}- \theta \delta.
\end{align*}
But by Proposition \ref{strongintB}, 
$$\log \EE  e^{\sup_{ \xi \in  W}  \{ \langle \theta \xi , X\rangle - \Lambda_p(\theta\xi) \}} \leq \kappa b(\theta W),$$
where $\kappa$ is a numerical constant,  and $b(\theta W)$ is defined in \eqref{Rademwidth}. Since $b(\theta W) = \theta b(W)$, we finally get
\begin{align*}
\log \PP\Big( \sup_{ \xi \in  W}  \{ \langle \theta \xi , X\rangle - \Lambda_p(\theta\xi) - \theta \delta\}& \geq \phi_p(t-\delta)-1\Big) \leq  - \phi_p(t-\delta) +1 \\
&+ \theta ( \kappa b(W)-  \delta).
\end{align*}
Thus, if $b(W) \leq \delta/\kappa$, we obtain
$$ \PP \big( f(X) \geq t \big) \leq |\mathcal{D}| e^{-\phi_p(t-\delta) +1}.$$
To complete the proof, it suffices to observe that the Rademacher mean-width of a set is the same as the one of its convex hull, so that $b(W) = b(V)$.

\section{Proof of Theorem \ref{revlogSob}}
 We write $\Lambda$ as a short-hand for $\Lambda_\mu$.
We will follow the lines of the argument in the Gaussian case from \cite [proof of Theorem 1]{LedouxE} which was based on Talagrand's transportation-entropy inequality \cite{Talagrand2}.  By definition,
$$\mathcal{I}(\nu) = \int \big( \langle \nabla \Lambda(\nabla f(x)), \nabla f(x) \rangle - \Lambda(\nabla f(x)) \big) e^{f(x)} d\mu(x).$$
Recall that we denote $\mu_{1/2} = \frac{1}{2} \delta_1 + \frac{1}{2} \delta_{-1}$. Let $i\in \{1,\ldots,n\}$. We have
$$\int \Lambda'(\partial_i f(x)) \partial_i f(x) e^{f(x)} d\mu_{1/2}(x_i) =\Lambda'(\partial_i f(x)) \partial_i f(x) \Big( \frac{1}{2}e^{f(x_+)} +\frac{1}{2}e^{f(x_-)} \Big),$$
where $x_+ = (x_1,\ldots,x_{i-1},1,x_{i-1},\ldots,x_n)$ and $x_- = (x_1,\ldots,x_{i-1},-1,x_{i+1},\ldots,x_n)$. But $\Lambda' = \tanh$, therefore
$$ \Lambda'(\partial_i f(x))= \frac{e^{f(x_+)} - e^{f(x_-)} }{e^{f(x_+)} + e^{f(x_-)}}.$$
Therefore,
\begin{align*}\int \Lambda'(\partial_i f(x)) \partial_i f(x) e^{f(x)} d\mu_{1/2}(x_i) &=\frac{1}{2}\big(e^{f(x_+)} - e^{f(x_-)} \big) \\
& =\int x_i  e^{f(x)} d\mu_{1/2}(x_i).
\end{align*}
Thus, 
$$\int \Lambda'(\partial_i f(x)) \partial_i f(x) e^{f(x)} d\mu_{1/2}(x_i) =\int x_i \partial_i f(x) e^{f(x)} d\mu_{1/2}(x_i).$$
Integrating the above equality with respect to $(x_1,\ldots,x_{i-1},x_{i+1},\ldots,x_n)$, and summing over $i\in\{1,\ldots,n\}$, we deduce 
$$ \int \langle\nabla \Lambda(\nabla f(x)), \nabla f(x) \rangle d\mu(x) = \int \langle x, \nabla f(x)\rangle d\mu(x).$$
Therefore,
$$ \mathcal{I}(\nu) = \int \big( \langle x, \nabla f(x)\rangle -\Lambda(\nabla f(x)) \big) d\mu(x).$$
In particular,
$$ \mathcal{I}(\nu) \leq  \int \sup_{\xi \in V} \big\{ \langle x, \xi \rangle -\Lambda(\xi) \big\} d\mu(x).$$
But, the Gibbs variational principle \eqref{Gibbs} implies that
$$  \int \sup_{\xi \in V} \big\{ \langle x, \xi \rangle -\Lambda(\xi) \big\} d\mu(x) \leq H(\nu|\mu) + \log \int e^{\sup_{\xi \in V} \{ \langle \xi, x\rangle - \Lambda(\xi)\}} d\mu(x).$$
Using Proposition \ref{strongintB}, we can conclude the proof.

\section*{Acknowledgment}
I thank Ofer Zeitouni for many fruitful discussions which helped me build the present paper, as well as his valuable comments on an earlier version of the manuscript. I am grateful to Ronen Eldan for several influential and helpful discussions.

\bibliographystyle{plain}
\bibliography{main}{}

\end{document}